# The MINI mixed finite element for the Stokes problem: An experimental investigation

Andrea Cioncolini[a]; Daniele Boffi[b,c]

[a] School of Mechanical, Aerospace and Civil Engineering, University of Manchester, George Begg Building, Sackville Street, M1 3BB Manchester, United Kingdom. E-mail: andrea.cioncolini@manchester.ac.uk

[b] Dipartimento di Matematica ''F. Casorati'', University of Pavia, Via Ferrata 1, 27100 Pavia, Italy.

[c] Department of Mathematics and System Analysis, Aalto University, Espoo, Finland.

E-mail: daniele.boffi@unipv.it

$O(h^{3/2})$ superconvergence in pressure and velocity has been experimentally investigated for the two-dimensional Stokes problem discretized with the MINI mixed finite element. Even though the classic mixed finite element theory for the MINI element guarantees linear convergence for the total error, recent theoretical results indicate that superconvergence of order $O(h^{3/2})$ in pressure and of the linear part of the computed velocity to the piecewise linear nodal interpolation of the exact velocity is in fact possible with structured, three-directional triangular meshes. The numerical experiments presented here suggest a more general validity of $O(h^{3/2})$ superconvergence, possibly to automatically generated and unstructured triangulations. In addition, the approximating properties of the complete computed velocity have been compared with the approximating properties of the piecewise-linear part of the computed velocity, finding that the former is generally closer to the exact velocity, whereas the latter conserves mass better.

**Keywords**: Stokes problem; Mixed finite element method; MINI finite element; Superconvergence; Benchmark numerical experiments.

## 1. Introduction

The Stokes problem represents the asymptotic limiting form of the Navier-Stokes problem for fluid dynamics when the Reynolds number becomes very small [1]. In this limit, the fluid dynamics is essentially controlled by diffusion, so that the non-linear convection term in the full Navier-Stokes equation can be dropped and the equation correspondingly linearized. Stokes flows, also called creeping flows, are typically characterized by small flow velocities, large fluid viscosities or small length scales, and are of relevance in lubrication theory, in porous media flow, in certain biological applications such as the swimming of microorganisms, and in microfluidics applications where fluid flows are geometrically constrained in the sub-millimeter scale. Besides their practical relevance, analytical solutions are available for several Stokes flows, thus making these problems essential and extensively used benchmark test cases in numerical fluid dynamics.

The present study is concerned with the numerical approximation of the two-dimensional Stokes problem using the MINI mixed finite element, which will be referred to as *Stokes-MINI* in the following. The analysis of mixed finite element methods for the Stokes problem can be based on the general theory of saddle point problems, developed by Babuška [2] and Brezzi [3]. The MINI mixed finite element, in particular, was introduced by Arnold, Brezzi and Fortin [4] specifically for the discretization of the Stokes problem and relies on continuous, piecewise-linear polynomials enriched with cubic bubble functions for the discrete velocity and on continuous, piecewise-linear polynomials for the discrete pressure. Standard mixed finite element theory [5] assures that Stokes-MINI is stable and converges linearly for both velocity and pressure, i.e. $O(h)$ convergence where $h$ is the discretization parameter. Despite the ''unbalanced'' approximation properties of the involved finite element spaces (the pressure space would allow for second order of approximation, while the velocity space is only linearly convergent), the Stokes-MINI element is quite popular due to its simplicity.

Although the standard error analysis guarantees linear convergence for Stokes-MINI, $O(h^{3/2})$ superconvergence in pressure and of the linear part of the computed velocity to the linear nodal interpolant of the exact velocity has been recently proved by Eichel, Tobiska and Xie [6] on structured, three-directional triangular meshes. Automatically generated, unstructured triangular meshes are not covered by the existing superconvergence theory, and this is what motivated the present work. The main objective of our study was to systematically investigate $O(h^{3/2})$ superconvergence for Stokes-MINI on a selection of seven benchmark test cases with analytical solution, using automatically generated unstructured triangular meshes. The test cases selected for use here are widely known in the literature and include four enclosed flows, the lid-driven cavity flow, and two more general flow problems with open boundary, inflow and outflow. Our results suggest a validity of $O(h^{3/2})$ superconvergence more general than what the existing theory covers, possibly to automatically generated unstructured triangular meshes. Contrary to common belief, our results also show that the bubble function not only stabilize the MINI finite element formulation, but also generally improves the quality of the velocity approximation, though at the expense of worsening the local mass conservation.

The rest of this paper is organized as follows: Section 2 includes the problem description and provides the necessary background material, Section 3 presents the Stokes flow benchmark test cases selected for use here, Section 4 describes the numerical procedure, while the results of the simulations are presented and discussed in Section 5. Particular care has been devoted to the description of the implementation of the test cases in the sake of reproducibility of the presented results.

## 2. Problem description

### 2.1. Strong formulation of the Stokes problem

The Stokes problem formulation of interest here reads as follows: find a two-dimensional velocity vector field $\underline{u}(x,y) = \big(u_x(x,y), u_y(x,y)\big)$, with $u_x, u_y \in C^2(\Omega) \cap C^0(\overline{\Omega})$ and a two-dimensional scalar pressure field $P(x,y) \in C^1(\Omega)$ such that:

$$-\mu \Delta \underline{u} + \underline{\nabla} P = \rho \underline{f} \quad in\ \Omega \tag{1}$$

$$div\ (\underline{u}) = 0 \quad in\ \Omega \tag{2}$$

$$\underline{u} = \underline{0} \quad on\ \partial\Omega \tag{3}$$

where $\Omega \subset \mathbb{R}^2$ is a bounded and connected polygonal domain in the plane with boundary $\partial\Omega$, $\rho$ and $\mu$ are the fluid density and viscosity (both assumed constant), while $\underline{f} = (f_x, f_y)$ with $f_x, f_y \in C^0(\Omega)$ is the external force field density per unit fluid mass. Problem (1)-(2) is appropriate for describing steady-state, isothermal and incompressible creeping flow of Newtonian fluids with constant density and constant viscosity, and can be easily generalized to density and viscosity given functions of position. As shown in Eq. (3), we consider a homogeneous Dirichlet boundary condition for the velocity. As is well known [5], a Stokes problem with a non-homogeneous Dirichlet boundary condition for the velocity can be restated as a homogeneous problem with a change of variable and modification to the right hand side of Eqs. (1) and (2), so that the problem formulation considered here also covers the non-homogeneous case. The Dirichlet boundary condition for the velocity enables modeling the interaction of the fluid with solid boundaries and the specification of inflows and outflows, thus allowing the treatment of a good deal of practical problems in fluid dynamics. It is well known [1,5,7] that in incompressible flows with the velocity specified everywhere on the boundary the pressure field is only determined up to an arbitrary additive constant. This implies that if a solution $(u_x, u_y, P)$ of the Stokes problem (1)-(3) exists then it is not unique, because $(u_x, u_y, P + c)$ will also be a solution for any constant $c \in \mathbb{R}$. Enforcing a null mean value of the pressure field over the entire domain $\Omega$ restores uniqueness:

$$\frac{1}{|\Omega|} \int_\Omega P \, d\Omega = 0 \qquad (4)$$

The unique solvability of the Stokes problem, and more generally of the Navier-Stokes problem with general boundary conditions, has not been proved yet and is still the subject of current investigation.

### 2.2. Weak formulation of the Stokes problem

The weak formulation of the Stokes problem (1)-(4) is derived by multiplying the Eqs. (1) and (2) by test functions $\underline{v} = (v_x, v_y)$; $q$ and then integrating on the domain $\Omega$ [5,7]. Using the Green's First Identity and the formula of integration by parts to transfer part of the derivatives onto the test functions yields the weak formulation of the Stokes problem: find a velocity vector field $\underline{u}(x, y) = \big(u_x(x, y), u_y(x, y)\big)$, with $u_x, u_y \in H_0^1(\Omega)$ and a scalar pressure field $P(x, y) \in L_0^2(\Omega)$ such that:

$$a(\underline{u}, \underline{v}) + b(\underline{v}, P) = F(\underline{v}) \quad \forall \, \underline{v} \in H_0^1(\Omega)^2 \qquad (5)$$

$$b(\underline{u}, q) = 0 \quad \forall \, q \in L_0^2(\Omega) \qquad (6)$$

where the forms $a(*,*), b(*,*)$ and the functional $F(*)$ are defined as follows:

$$a(\underline{u}, \underline{v}) = \mu \int_\Omega \underline{\nabla}\underline{u} : \underline{\nabla}\underline{v} \, d\Omega = \mu \int_\Omega \underline{\nabla} u_x \cdot \underline{\nabla} v_x \, d\Omega + \mu \int_\Omega \underline{\nabla} u_y \cdot \underline{\nabla} v_y \, d\Omega \qquad (7)$$

$$b(\underline{v}, P) = -\int_\Omega div(\underline{v}) P \, d\Omega \qquad (8)$$

$$F(\underline{v}) = \rho \int_\Omega \underline{f} \cdot \underline{v} \, d\Omega \qquad (9)$$

where $\underline{f} = (f_x, f_y)$ with $f_x, f_y \in L^2(\Omega)$ is the external force field density per unit fluid mass. Here $H_0^1(\Omega)$ is the usual Sobolev space of functions that are square-integrable with their first weak derivatives and that vanish at the domain boundary in the sense of traces, while $L_0^2(\Omega)$ is the subspace of $L^2(\Omega)$ of square-integrable functions with vanishing mean. As is well known [5], the weak Stokes problem (5)-(6) admits a unique weak solution $(\underline{u}, P) \in H_0^1(\Omega)^2 \times L_0^2(\Omega)$.

### 2.3. Discrete Galerkin approximation of the Stokes problem

The discrete Galerkin approximation of the weak Stokes problem (5)-(6) consists of formulating the problem in two families of finite-dimensional linear subspaces $X_h(\Omega) \subset H_0^1(\Omega)$ and $M_h(\Omega) \subset L_0^2(\Omega)$ that depend on a discretization parameter $h$ and approximate the Hilbert spaces $H_0^1(\Omega)$ and $L_0^2(\Omega)$, so that the discrete Stokes problem reads as follows: find $(\underline{u}_h, P_h) \in X_h^2 \times M_h$ such that:

$$a(\underline{u}_h, \underline{v}_h) + b(\underline{v}_h, P_h) = F(\underline{v}_h) \quad \forall \, \underline{v}_h \in X_h(\Omega)^2 \tag{10}$$

$$b(\underline{u}_h, q_h) = 0 \quad \forall \, q_h \in M_h(\Omega) \tag{11}$$

The discrete problem (10)-(11) admits a unique solution provided that the discrete subspaces $X_h$ and $M_h$ satisfy the following discrete inf-sup condition [5]:

$$\exists \, \beta > 0: \quad \inf_{0 \neq q_h \in M_h} \sup_{0 \neq \underline{v}_h \in X_h^2} \frac{b(\underline{v}_h, q_h)}{\| \underline{v}_h \|_X \, \| q_h \|_M} \geq \beta \tag{12}$$

where $\beta$ is independent of the discretization parameter $h$. The discrete inf-sup condition (12) acts as a compatibility condition between the spaces $X_h$ and $M_h$, and can therefore be regarded as a recipe for the construction of well-posed schemes that assure unique solvability of the discrete Stokes problem.

*2.4. The MINI mixed finite element for the Stokes problem*

When approximating the weak Stokes problem with mixed finite elements, the finite-dimensional subspaces $X_h$ and $M_h$ consist of piecewise polynomials constructed on a discretization $\Omega_h$ of the domain $\Omega$. Here, in particular, the attention is restricted to triangulations $\mathcal{T}_h$ of the domain $\Omega$ that are conformal (two triangles share at most one vertex or one edge) and regular (the smallest angle in all triangles is bounded away from zero independently of the mesh parameter). As already noted, the domain $\Omega$ is here assumed bounded, connected and polygonal: it is therefore possible to divide its closure $\overline{\Omega}$ into triangles $T$ that form a triangulation $\mathcal{T}_h$ that wholly covers $\Omega$, so that $\Omega_h = \mathcal{T}_h \equiv \Omega$, and where the positive mesh spacing parameter $h$ represents the longest edge across all triangles in the triangulation.

Several pairs of mixed finite elements have been proposed for the Stokes problem, and more generally for numerical fluid dynamics applications, where different mixed finite element schemes basically differ in the global regularity and local order of the polynomials. It is common practice to schematically identify mixed finite element pairs with the symbol $[\mathbb{P}_k]^d/\mathbb{P}_m$, which specifies the number of space dimensions $d$ (2 in the present case) and the local order of the polynomials used for the velocity ($k$) and pressure ($m$) spaces. The MINI finite element was introduced by Arnold, Brezzi and Fortin [4] specifically for the discretization of the Stokes problem and relies on continuous, piecewise-linear polynomials enriched with local cubic bubble functions for the discrete velocity space $X_h$ and on continuous, piecewise-linear polynomials for the discrete pressure space $M_h$, so that it is characterized by the $[\mathbb{P}_{1+b}]^2/\mathbb{P}_1$ pair. As is well known [5], the $[\mathbb{P}_1]^2/\mathbb{P}_1$ pair that would be computationally quite convenient is however unstable, because it fails to satisfy the discrete inf-sup condition (12). The idea behind the MINI element was therefore to stabilize the unstable $[\mathbb{P}_1]^2/\mathbb{P}_1$ pair by adding local functions, named *bubbles*, to properly enrich the discrete velocity space. In the MINI element, in particular, the bubble function is a cubic polynomial defined locally in each triangle and given, as shown in Eq. (14), by the product of the linear nodal basis functions (baricentric coordinates) $\varphi_1, \varphi_2, \varphi_3$ of the triangle itself:

$$X_h = \{v_h \in C^0(\bar{\Omega}): v_{h|T} \in \mathbb{P}_{1+b} \ \forall T \in \mathcal{T}_h\} \tag{13}$$

$$v_{h|T} = a_h + b_h x + c_h y + d_h\, \varphi_1(x,y)\varphi_2(x,y)\varphi_3(x,y) \tag{14}$$

$$M_h = \{q_h \in C^0(\bar{\Omega}): q_{h|T} \in \mathbb{P}_1 \ \forall T \in \mathcal{T}_h\} \tag{15}$$

$$q_{h|T} = a_h + b_h x + c_h y \tag{16}$$

where $a_h, b_h, c_h, d_h$ are constants while $v_{h|T}$ in Eq. (14) refers to one of the discrete velocity components. Since with the MINI finite element the velocity is approximated with continuous, piecewise-linear polynomials enriched with bubble functions, it is possible to split the computed velocity $\underline{u}_h$ as follows:

$$\underline{u}_h = \underline{u}_{hl} + \underline{u}_{hb} \tag{17}$$

where $\underline{u}_{hl}$ is the piecewise-linear part of the computed velocity while $\underline{u}_{hb}$ is the bubble part. As originally pointed out by Verfürth [8] and successively reiterated by Bank and Welfert [9,10], Kim and Lee [11] and Russo [12], the piecewise-linear part of the computed velocity $\underline{u}_{hl}$ actually seems to be a better approximation to the exact velocity $\underline{u}$ than the complete computed velocity $\underline{u}_h$ itself. The bubble part $\underline{u}_{hb}$, therefore, is only needed to stabilize the formulation but apparently does not improve the quality of the velocity approximation. This is the reason why several a-posteriori error estimators for the Stokes-MINI problem are actually based on the linear part of the computed velocity [8,12].

Standard mixed finite element theory [5] assures that Stokes-MINI is stable and converges linearly:

$$\| \underline{u} - \underline{u}_h \|_{H^1} + \| P - P_h \|_{L^2} \leq C\, h\, \left( \| \underline{u} \|_{H^2} + \| P \|_{H^1} \right) \tag{18}$$

where $C$ is a positive constant independent of $h$, $\|*\|_{H^1}$ and $\|*\|_{L^2}$ are the usual norms in $H^1(\Omega)$ and $L^2(\Omega)$, and provided the exact solution $(\underline{u}, P) \in H^2(\Omega)^2 \times H^1(\Omega)$. Even though the standard error analysis assures linear convergence for the global error, Eichel, Tobiska and Xie [6] recently proved the following $O(h^{3/2})$ superconvergence result on three-directional structured triangular meshes (i.e. triangular meshes generated starting from a structured rectangular mesh and dividing each rectangle into two triangles using one of the rectangle diagonals):

**THEOREM 1.** *With reference to the Stokes-MINI problem, assume that the triangulation $\mathcal{T}_h$ is three-directional structured, and that the exact solution $(\underline{u}, P) \in H^3(\Omega)^2 \times H^2(\Omega)$; then:*

$$\left( |\underline{u}_{hl} - i_h \underline{u}|_{H^1}^2 + \| P_h - j_h P \|_{L^2}^2 \right)^{1/2} \leq C h^{3/2} \left( \| \underline{u} \|_{H^3} + \| P \|_{H^2} \right) \tag{19}$$

where $C$ is a positive constant independent of $h$, $| * |_{H^1}$ is the usual semi-norm in $H^1(\Omega)$, and $(i_h \underline{u}, j_h P)$ denotes the standard, piecewise-linear nodal interpolation of the exact solution $(\underline{u}, P)$.

As can be seen, Theorem 1 guarantees superconvergence of order $O(h^{3/2})$ of the liner part of the computed velocity $\underline{u}_{hl}$ to the piecewise-linear nodal interpolation of the exact velocity $i_h \underline{u}$, and of the computed pressure $P_h$ to the piecewise-linear nodal interpolation of the exact pressure $j_h P$. As noted by the authors [6], the superconvergence in pressure can be generalized as follows:

$$\| P - P_h \|_{L^2} \leq \| P_h - j_h P \|_{L^2} + \| j_h P - P \|_{L^2} \leq C \left( h^{3/2} + h^2 \right) \left( \| \underline{u} \|_{H^3} + \| P \|_{H^2} \right) \tag{20}$$

which, under the same assumptions of Theorem 1, guarantees superconvergence of order $O(h^{3/2})$ of the computed pressure $P_h$ to the exact pressure $P$. It is worth noting that any superconvergence in pressure is not in conflict with the general error bound (18): the velocity error does converge linearly, and this controls the rate of convergence of the global error hiding any superconvergence in pressure, unless the pressure error is considered alone.

As previously noted, the weak formulation of the Stokes problem of interest here is set in the Sobolev space $H_0^1(\Omega)$. As is well known, in $H_0^1(\Omega)$ the semi-norm $| * |_{H^1}$ is itself a norm, and is equivalent to the usual norm $\| * \|_{H^1}$. This equivalence is exploited in the following: convergence rates in the sense of Theorem 1 of the linear part of the computed velocity $\underline{u}_{hl}$ to the piecewise linear nodal interpolation of the exact velocity $i_h \underline{u}$ are assessed using the norm $\| * \|_{H^1}$.

## 3. Benchmark test cases

The seven test cases with analytical solution selected for use here are well-known problems in computational fluid dynamics and are described below, while their corresponding velocity vector fields are provided in Fig. 1. Notably, four test cases (#1, #2, #3, #4) have homogeneous Dirichlet boundary condition for the velocity, while three (#5, #6, #7) have non-homogeneous boundary conditions. Moreover, four test cases (#1, #2, #5, #7) are polynomial, while three (#3, #4, #6) are not. In all test problems the analytical pressure solution has vanishing mean over the domain.

### 3.1. Test problem #1 (enclosed vortex, polynomial)

Consider the Stokes problem (1)-(2) on the unit square domain $\Omega = (0,1)\times(0,1)$ with right-hand side in Eq. (1) given as follows:

$$\rho f_x = -\mu \{4y(1-y)(2y-1)[(1-2x)^2 - 2x(1-x)] + 12x^2(1-x)^2(1-2y)\} \tag{21}$$
$$+(1-2x)(1-y)$$
$$\rho f_y = -\mu \{4x(1-x)(1-2x)[(1-2y)^2 - 2y(1-y)] + 12y^2(1-y)^2(2x-1)\} \tag{22}$$
$$-x(1-x)$$

and with homogeneous Dirichlet boundary condition for the velocity on all domain boundary:

$$u_x = u_y = 0 \text{ on } \partial\Omega \tag{23}$$

The corresponding exact solution is:

$$u_x(x,y) = x^2(1-x)^2 2y(1-y)(2y-1) \tag{24}$$
$$u_y(x,y) = y^2(1-y)^2 2x(1-x)(1-2x) \tag{25}$$
$$P(x,y) = x(1-x)(1-y) - 1/12 \tag{26}$$

As shown in Fig. 1, this flow problem corresponds to an enclosed flow where the velocity field has the form of a large vortex rotating clock-wise.

### 3.2. Test problem #2 (enclosed vortex, polynomial)

Consider the Stokes problem (1)-(2) on the unit square domain $\Omega = (0,1)\times(0,1)$ with right-hand side in Eq. (1) given as follows:

$$\rho f_x = -\mu \{(2 - 12x + 12x^2)(2y - 6y^2 + 4y^3) + (x^2 - 2x^3 + x^4)(-12 + 24y)\} + 1/24 \tag{27}$$
$$\rho f_y = \mu \{(2 - 12y + 12y^2)(2x - 6x^2 + 4x^3) + (y^2 - 2y^3 + y^4)(-12 + 24x)\} + 1/24 \tag{28}$$

and with homogeneous Dirichlet boundary condition for the velocity on all domain boundary:

$$u_x = u_y = 0 \text{ on } \partial\Omega \tag{29}$$

The corresponding exact solution is:

$$u_x(x,y) = (x^2 - 2x^3 + x^4)(2y - 6y^2 + 4y^3) \tag{30}$$
$$u_y(x,y) = -(2x - 6x^2 + 4x^3)(y^2 - 2y^3 + y^4) \tag{31}$$
$$P(x,y) = (x + y - 1)/24 \tag{32}$$

As shown in Fig. 1, this flow problem corresponds to an enclosed flow where the velocity field has the form of a large vortex rotating counter clock-wise.

### 3.3. Test problem #3 (enclosed vortex, non-polynomial)

Consider the Stokes problem (1)-(2) on the unit square domain $\Omega = (0,1)\times(0,1)$ with right-hand side in Eq. (1) given as follows:

$$\rho f_x = -4\pi^2 \mu \sin(2\pi y)[2\cos(2\pi x) - 1] + 4\pi^2 \sin(2\pi x) \tag{33}$$

$$\rho f_y = 4\pi^2 \mu \sin(2\pi x)[2\cos(2\pi y) - 1] - 4\pi^2 \sin(2\pi y) \tag{34}$$

and with homogeneous Dirichlet boundary condition for the velocity on all domain boundary:

$$u_x = u_y = 0 \text{ on } \partial\Omega \tag{35}$$

The corresponding exact solution is:

$$u_x(x,y) = \sin(2\pi y)[1 - \cos(2\pi x)] \tag{36}$$

$$u_y(x,y) = \sin(2\pi x)[\cos(2\pi y) - 1] \tag{37}$$

$$P(x,y) = 2\pi[\cos(2\pi y) - \cos(2\pi x)] \tag{38}$$

As shown in Fig. 1, this flow problem corresponds to an enclosed flow where the velocity field has the form of a large vortex rotating counter clock-wise.

### 3.4. Test problem #4 (enclosed vortex, non-polynomial)

Consider the Stokes problem (1)-(2) on the unit square domain $\Omega = (0,1)\times(0,1)$ with right-hand side in Eq. (1) given as follows:

$$\begin{aligned}\rho f_x = &-\mu\{2e^x[(x^2 + x - 1)(x^2 + 3x - 2) + (x^2 - x)(2x + 3)](y^2 - y)(2y - 1) \\ &+ 2e^x(x-1)^2 x^2(12y - 6)\} \\ &+ (y^2 - y)e^x[x^4(y^2 - y + 12) + 6x^3(y^2 - y - 4) + x^2(y^2 - y + 12) \\ &+ 8x(y - y^2) + 2y^2 - 2y]\end{aligned} \tag{39}$$

$$\begin{aligned}\rho f_y = &-\mu\{-e^x(x^4 + 10x^3 + 19x^2 - 6x - 6)(y-1)^2 y^2 \\ &- 2e^x(x^2 - x)(x^2 + 3x - 2)(6y^2 - 6y + 1)\} \\ &+ (2y \\ &- 1)\{-456 \\ &+ e^x[x^4(y^2 - y + 12) + 2x^3(y^2 - y - 36) + x^2(-5y^2 + 5y + 228) \\ &+ 2x(y^2 - y - 228) + 456]\} + e^x(x^4 + 2x^3 - 5x^2 + 2x)(2y - 1)(y^2 - y)\end{aligned} \tag{40}$$

And with homogeneous Dirichlet boundary condition for the velocity on all domain boundary:

$$u_x = u_y = 0 \text{ on } \partial\Omega \tag{41}$$

The corresponding exact solution is:

$$u_x(x, y) = 2e^x(x - 1)^2 x^2 (y^2 - y)(2y - 1) \tag{42}$$

$$u_y(x, y) = -e^x(x^2 - x)x(x^2 + 3x - 2)(y - 1)^2 y^2 \tag{43}$$

$$P(x, y) = -424 + 156e \tag{44}$$

$$+ (y^2 - y)\{-456$$

$$+ e^x[x^4(y^2 - y + 12) + 2x^3(y^2 - y - 36) + x^2(-5y^2 + 5y + 228)$$

$$+ 2x(y^2 - y - 228) + 456]\}$$

As shown in Fig. 1, this flow problem corresponds to an enclosed flow where the velocity field has the form of a large and slightly asymmetric vortex rotating counter clock-wise.

### 3.5. Test problem #5 (lid-driven cavity flow, polynomial)

Consider the Stokes problem (1)-(2) on the unit square domain $\Omega = (0,1) \times (0,1)$ with right-hand side in Eq. (1) given as follows:

$$\rho f_x = 0 \tag{45}$$

$$\rho f_y = \mu[(12x - 6)(y^4 - y^2) + (8x^3 - 12x^2 + 4x)(6y^2 - 1) + 0.4(6x^5 - 15x^4 + 10x^3)] \tag{46}$$

The boundary conditions for the velocity are of the homogeneous Dirichlet type on all domain boundaries except along the top edge of the domain, where the velocity is given as indicated below:

$$u_x(x, 1) = x^4 - 2x^3 + x^2; \; u_y(x, 1) = 0 \tag{47}$$

The corresponding exact solution is:

$$u_x(x, y) = (x^4 - 2x^3 + x^2)(2y^3 - y) \tag{48}$$

$$u_y(x, y) = -(2x^3 - 3x^2 + x)(y^4 - y^2) \tag{49}$$

$$P(x, y) = \mu[(4x^3 - 6x^2 + 2x)(2y^3 - y) + 0.4(6x^5 - 15x^4 + 10x^3)y - 0.1] \tag{50}$$

As shown in Fig. 1, this flow problem corresponds to a lid-driven cavity flow where the velocity field is driven by a specified body force in addition to the non-uniform shear acting on the top edge of the domain boundary. This flow problem is a variant of the lid-driven cavity flow, a well-known benchmark validation case for numerical fluid dynamics, proposed by Shih et al. [13]. Differently from other versions of the lid-driven cavity flow, this flow problem has analytical solution and no velocity singularities at the top corners of the domain.

### 3.6. Test problem #6 (corner flow, non-polynomial)

Consider the Stokes problem (1)-(2) on the unit square domain $\Omega = (0,1) \times (0,1)$ with right-hand side in Eq. (1) given as follows:

$$\rho f_x = -\mu[sin(xy)x(x^2 + y^2) - 2cos(xy)y] - sin(xy)y \tag{51}$$

$$\rho f_y = \mu[sin(xy)y(x^2 + y^2) - 2cos(xy)x] - sin(xy)x \tag{52}$$

The boundary conditions for the velocity are of the homogeneous Dirichlet type on all domain boundaries except along the top and right edges of the domain, where the velocity is given as indicated below:

$$\begin{aligned} u_x(x,1) = -xsin(x); \; u_y(x,1) = sin(x) \\ u_x(1,y) = -sin(y); \; u_y(1,y) = ysin(y) \end{aligned} \tag{53}$$

The corresponding exact solution is:

$$u_x(x,y) = -sin(xy)x \tag{54}$$

$$u_y(x,y) = sin(xy)y \tag{55}$$

$$P(x,y) = cos(xy) - 0.9460830703671845 \tag{56}$$

As shown in Fig. 1, this flow problem corresponds to a corner flow entering the domain from the right edge and exiting the domain through the top edge.

### 3.7. Test problem #7 (colliding flow, polynomial)

Consider the Stokes problem (1)-(2) on the square domain $\Omega = (-1,1)\times(-1,1)$ with zero body force:

$$\rho f_x = \rho f_y = 0 \tag{57}$$

The boundary conditions for the velocity are of Dirichlet type as indicated below:

$$\begin{aligned} u_x(x,1) = 20x - 4x^5; \; u_y(x,1) = 20x^4 - 4 \\ u_x(1,y) = 20y^4 - 4; \; u_y(1,y) = 20y - 4y^5 \\ u_x(x,-1) = 20x - 4x^5; \; u_y(x,-1) = -20x^4 + 4 \\ u_x(-1,y) = -20y^4 + 4; \; u_y(-1,y) = 20y - 4y^5 \end{aligned} \tag{58}$$

The corresponding exact solution is:

$$u_x(x,y) = 20xy^4 - 4x^5 \tag{59}$$

$$u_y(x,y) = 20x^4y - 4y^5 \tag{60}$$

$$P(x,y) = \mu(120x^2y^2 - 20x^4 - 20y^4 - 32/6) \tag{61}$$

As shown in Fig. 1, this flow problem corresponds to a colliding flow entering the domain from the midpoints of all edges and exiting the domain through the domain corners.

## 4. Numerical methods

All calculations have been performed with purpose-written scripts in MATLAB R2017a (64-bit), with the kinematic viscosity $\mu$ in Eq. (1) set to 1 (kg/ms). Details on the mesh generation, numerical quadrature and linear system solution are provided below.

### 4.1. Mesh generation

The mesh generation was performed using DistMesh, a simple open-source MATLAB code for the generation of unstructured triangular (and tetrahedral) meshes developed by Persson and Strang [14]. As a representative example, a uniformly spaced mesh generated in the unit square using DistMesh is provided in Fig. 2 (left). All meshes used in the present work were generated with uniform spacing. Sometimes, DistMesh produces meshes with corner triangles, which are triangles that have two edges on the domain boundary. In the example in Fig. 2 (left), in particular, this happens in the top-right and bottom-right corners. As is well known, corner triangles are undesirable in low-order finite element settings because they are strongly influenced by the boundary condition, and should therefore be avoided. In the numerical tests presented here, therefore, corner triangles were eliminated via diagonal exchange, as illustrated in Fig. 2 (right). DistMesh implements an iterative technique based on a physical analogy between a simplex mesh and a truss structure that typically produces triangles that are almost equilateral, and therefore meshes with particularly good quality, thus avoiding deformed triangles with large or small angles that adversely affect error estimates. This is, in fact, what motivated the use of DistMesh in the present work. Here, in particular, the mesh quality was double-checked using the following widely used triangle quality measures [15]:

$$q_1 = \frac{(b+c-a)(c+a-b)(a+b-c)}{abc} \tag{62}$$

$$q_2 = \frac{4\sqrt{3}\,A_T}{a^2+b^2+c^2} \tag{63}$$

where $a, b, c$ are the triangle side lengths and $A_T$ is the triangle area. An equilateral triangle has $q_1 = q_2 = 1$ while a degenerate, zero area triangle has $q_1 = q_2 = 0$. As a rule-of-thumb, in a good quality mesh all triangles should have $q_1, q_2$ above about 0.4-0.5. All meshes generated in the present study had quality measures $q_1, q_2$ always above 0.7 and with most triangles actually above 0.8-0.9, and can therefore be regarded as good quality meshes. As a representative example, Fig. 3 provides the quality histograms for the mesh in the unit square presented in Fig. 2 (right): as can be seen, the quality of the mesh is quite good.

As is well known, Stokes flows in confined domains present vortices, known as *Moffatt vortices*, near corners between intersecting solid boundaries or between a solid boundary and a free surface. Very fine grids are required to capture Moffatt vortices in numerical simulations [16], way finer than the meshes used in the present work. As such, the occurrence of Moffatt vortices is not investigated in the test problems presented in what follows.

*4.2. Numerical integration*

The evaluation of the linear system matrix and vector components, as well as the evaluation of the errors, requires the numerical integration of various functions over the mesh triangles. Depending on the problem being considered, these integrals might involve polynomial or non-polynomial integrand functions. As is well known, the accurate numerical evaluation of such integrals is an essential prerequisite in any finite element discretization, so that special care was devoted to the selection and implementation of appropriate quadrature formulas. Among the several numerical integration formulas available in the literature, the quadrature formulas developed by Dunavant [17] were chosen for use here. These are symmetrical Gaussian quadrature formulas of degree up to 20, specifically developed for the numerical integration over triangles for use in finite element formulations. In the numerical tests presented in the following, the quadrature formula was always selected of appropriate degree for polynomial integrands (i.e. quadrature rule of degree no less than the degree of the polynomial integrand function), while the highest degree formula (degree 20) was always used for non-polynomial integrand functions.

*4.3. Linear system solution*

The linear system was solved using the Generalized Minimum Residual Method (GMRES) with incomplete LU factorization (ILU) as preconditioner, using MATLAB built-in functions. Several methods are available to solve linear systems. Among these, Krylov subspace iterative methods are regarded as one of the most efficient procedures currently available to solve symmetric and indefinite linear systems such as those of interest here [18]. As well known, the efficiency of Krylov subspace methods heavily relies on appropriate preconditioning strategies [19]. Among Krylov subspace methods, GMRES with ILU preconditioning can be regarded as a '*go-to*' technique that can effectively handle many practical problems [20], and this motivated its use here. The Minimum Residual Method (MINRES), which is quite popular for the Stokes problem and notably exploits the symmetry of the system matrix (which GMRES does not), requires however more effort to design an efficient preconditioner, and was therefore not considered here as a first option. In all calculations presented here, the initial guess was the zero vector (the default in MATLAB), while the drop tolerance in the ILU factorization was chosen by trial-and-error in the range $10^{-4}$-$10^{-1}$. Within the limits of the present work, GMRES with ILU preconditioning worked remarkably well. In particular, for the Test Problems #1 through #5 with no inflow or outflow convergence was achieved within a few tens of iterations, and the returned iteration had relative residuals in the range of $10^{-15}$ to $10^{-10}$, which can be regarded as quite satisfactory. On the other hand, achieving convergence in Test Problems #6 and #7 with inflow/outflow required some fine-tuning of the drop tolerance and more iterations, and the returned iteration had relative residuals in the range of $10^{-7}$ to $10^{-5}$: still acceptable but clearly not as good as with the other test cases. This is no surprise, as incompressible flow problems with inflow/outflow are well known to be challenging to solve. Open boundaries must allow fluid to enter/leave the domain ensuring global mass conservation, which in incompressible flow problems implies that the amount of fluid mass (or volume) entering the domain must match exactly the mass (or volume) of fluid flowing out. Even when this is exactly satisfied in the continuous problem, it is only approximated in the discrete problem, and this makes incompressible flow problems with inflow/outflow boundary conditions particularly challenging to handle at the discrete level, and correspondingly more difficult to solve. As a matter of fact, incompressible flow problems with inflow/outflow boundary conditions are frequently reformulated using alternative boundary conditions of the Neumann type [7].

*4.4. Assessment*

The main objective of the present work was to experimentally ascertain $O(h^{3/2})$ superconvergence in the sense of Theorem 1 previously discussed on automatically generated, unstructured triangular meshes, using the seven benchmark test problems previously described. Automatically generated, unstructured triangular meshes are not covered by the existing superconvergence theory in Theorem 1, which is restricted to three-directional structured triangulations. Our experimental assessment of superconvergence in this more general setting will be instrumental in informing future theoretical developments of superconvergence for Stokes-MINI. Operatively, evaluating the order of convergence in the $H^1$ norm (equivalent in the present setting to the $H^1$ semi-norm, as previously noted) of the linear part of the computed velocity to the piecewise-linear nodal interpolation of the exact velocity, and the order of convergence in the $L^2$ norm of the computed pressure to the exact pressure, assessed superconvergence.

A further objective of the present work was to compare the approximating properties of the piecewise-linear part of the computed velocity $\underline{u}_{hl}$ to the approximating properties of the complete computed velocity $\underline{u}_h = \underline{u}_{hl} + \underline{u}_{hb}$, which also includes the bubble function $\underline{u}_{hb}$, again using the benchmark test problems previously described.

As previously noted, the piecewise-linear part of the computed solution $\underline{u}_{hl}$ is currently believed to be a better approximation to the exact solution $\underline{u}$ than the complete computed solution $\underline{u}_h$ itself [8-12], so that the bubble function is only needed to stabilize the formulation but apparently does not improve the quality of the velocity approximation. To the best of our knowledge, however, this result is presently not quantified in the open literature, and this motivated the second objective of the present work. Our experimental results in this respect will inform future a-posteriori error estimation studies for the Stokes-MINI problem, which following the common belief are currently based on the linear part of the computed velocity [8,12]. Operatively, the approximating properties of the linear part of the computed velocity and of the complete computed velocity were assessed in two steps. First, the respective approximating errors in the $H^1$ and $L^2$ norms were compared: $\|\underline{u} - \underline{u}_h\|_{H^1}$ was compared with $\|\underline{u} - \underline{u}_{hl}\|_{H^1}$, while $\|\underline{u} - \underline{u}_h\|_{L^2}$ was compared with $\|\underline{u} - \underline{u}_{hl}\|_{L^2}$. Second, the $L^2$ norms of the divergence of the computed velocities were compared, so that $\|div(\underline{u}_h)\|_{L^2}$ was compared with $\|div(\underline{u}_{hl})\|_{L^2}$. As formalized in Eq. (2), the exact velocity solution $\underline{u}$ is pointwise divergence-free. In the weak formulation, this property is equivalent to $\|div(\underline{u})\|_{L^2} = 0$. As several other low-order mixed finite elements, the MINI element is not weakly divergence-free, so that the $L^2$ norm of the divergence of the computed velocity is not zero. This quantity, in fact, converges to zero as the mesh is gradually refined, and its absolute value can be regarded as a measure of how closely the computed velocity conserves mass. Therefore, comparing $\|div(\underline{u}_h)\|_{L^2}$ with $\|div(\underline{u}_{hl})\|_{L^2}$ allows assessing which one, $\underline{u}_h$ or $\underline{u}_{hl}$, conserves mass better.

## 5. Results and discussion

Convergence histories for the velocity and pressure errors are provided in Fig. 4 for all test problems, while the corresponding convergence rates are summarized in Table 1.

As can be seen in Table 1 (second and third columns from the left), the $L^2$ and $H^1$ norms of the velocity error, $\|\underline{u} - \underline{u}_h\|_{L^2}$ and $\|\underline{u} - \underline{u}_h\|_{H^1}$, respectively converge with rates within 1.96-2.13 and 1.02-1.07 that appear consistent with the $O(h^2)$ and $O(h)$ convergence rates expected from mixed finite element theory [5], giving therefore confidence in the present numerical approach and implementation. The $H^1$ norm of the velocity error $\|i_h\underline{u} - \underline{u}_{hl}\|_{H^1}$ (third column from the right in Table 1) calculated with the piecewise-linear interpolation of the exact velocity and the linear part of the computed velocity converges with rates within 1.32-1.67 that appear consistent with $O(h^{3/2})$ convergence (to within ± 12%). Similarly for the $L^2$ norm of the pressure error $\|P - P_h\|_{L^2}$ (fourth column from the left in Table 1), which converges with rates within 1.41-1.59 that again appear consistent with $O(h^{3/2})$ convergence (to within ± 6%). Since the existing $O(h^{3/2})$ superconvergence

theory embodied in Theorem 1 previously discussed is only valid for three-directional structured triangular meshes, the present results suggest a validity of the $O(h^{3/2})$ superconvergence more general than what the existing theory covers, possibly to automatically generated unstructured triangular meshes of the type used here. Interestingly, the L² norm of the velocity error $\| i_h \underline{u} - \underline{u}_{hl} \|_{L^2}$ calculated with the piecewise-linear interpolation of the exact velocity and the linear part of the computed velocity (fourth column from the right in Table 1) converges with rates within 1.95-2.22 that appear consistent with $O(h^2)$ convergence. This indicates that $O(h^{3/2})$ superconvergence in velocity is only evident in the H¹ norm, while the rate of convergence in the L² norm is unaffected. Convergence rates for the velocity errors $\| \underline{u} - \underline{u}_{hl} \|_{L^2}$ and $\| \underline{u} - \underline{u}_{hl} \|_{H^1}$ are provided in Table 1 (last two columns on the right). As can be seen, the respective orders are within 1.95-2.12 and 1.00-1.04 for the L² and H¹ norms, and compare favorably with the corresponding convergence rates of $\| \underline{u} - \underline{u}_h \|_{L^2}$ and $\| \underline{u} - \underline{u}_h \|_{H^1}$ previously discussed (second and third columns from the left in Table 1), and also with the corresponding $O(h^2)$ and $O(h)$ convergence rates expected from mixed finite element theory [5].

The comparison between $\| \underline{u} - \underline{u}_h \|_{H^1}$ and $\| \underline{u} - \underline{u}_{hl} \|_{H^1}$ and between $\| \underline{u} - \underline{u}_h \|_{L^2}$ and $\| \underline{u} - \underline{u}_{hl} \|_{L^2}$ is presented in Fig. 5 (top), where the ratio of these errors is plotted as a function of the mesh spacing parameter, while the ratio of $\| div(\underline{u}_h) \|_{L^2}$ to $\| div(\underline{u}_{hl}) \|_{L^2}$ is presented as a function of the mesh spacing parameter in Fig. 5 (bottom). As can be seen in Fig. 5 (top-left), the velocity errors ratio in the H¹ norm gradually decreases as the mesh is gradually refined, with a trend that is similar for all test problems though the values of the velocity error ratio are clearly problem-specific. The values of the velocity errors ratio for the test problems with closed boundary (Test Problems #1 through #4), or with open boundary but no inflow/outflow (Test Problem #5), are strictly lower than one, indicating that for these problems the complete computed velocity $\underline{u}_h$ approximates the exact velocity $\underline{u}$ in the H¹ norm better than the piecewise-linear part of the computed velocity $\underline{u}_{hl}$, and in particular the better approximation of $\underline{u}_h$ with respect to $\underline{u}_{hl}$ becomes more evident as the mesh is gradually refined. Regarding the test problems with open boundary and inflow/outflow (Test Problems #6 and #7), the values of the velocity errors ratio indicate that for these problems the complete computed velocity $\underline{u}_h$ approximates the exact velocity $\underline{u}$ in the H¹ norm better than the piecewise-linear part of the computed velocity $\underline{u}_{hl}$ on fine meshes, while the converse is true on coarse meshes. As can be noticed in Fig. 5 (top-right), the velocity errors ratio in the L² norm are strictly lower than one, indicating that for all test problems the complete computed velocity $\underline{u}_h$ approximates the exact velocity $\underline{u}$ in the L² norm better than the piecewise-linear part of the computed velocity $\underline{u}_{hl}$. The dependence of the velocity errors ratio on the mesh spacing parameter is now

milder, and not detectable for Test Problem #6. Finally, as evident in Fig 5 (bottom), the L² norm of the divergence of the complete computed velocity ‖ $div(\underline{u}_h)$ ‖$_{L^2}$ is always bigger than that of the piecewise-linear part of the computed velocity ‖ $div(\underline{u}_{hl})$ ‖$_{L^2}$. This indicates that in all test problems the piecewise-linear part of the computed velocity $\underline{u}_{hl}$ conserves mass better than the complete computed velocity $\underline{u}_h$. In particular, the trend in the results in Fig. 5 (bottom) is decreasing as the mesh spacing parameter is gradually decreased, indicating that the improved mass conservation of $\underline{u}_{hl}$ with respect to $\underline{u}_h$ is more pronounced the coarser the mesh.

In conclusion, therefore, in all test problems the complete computed velocity $\underline{u}_h$ approximates the exact velocity in the L² norm better that the piecewise-linear part of the computed velocity $\underline{u}_{hl}$. This is also true in the H¹ norm, with the exception of flow problems with open boundary and inflow/outflow when solved on coarse meshes. The improved accuracy of the complete computed velocity, however, comes at the expense of a worse mass conservation. The piecewise-linear part of the computed velocity, in fact, conserves mass better than the complete computed velocity, particularly on coarse meshes.

The dependence of the results in Fig. 5 on the test problem can be traced back to the velocity boundary condition that characterize these problems, which in turn affected the accuracy of the numerical solution. As previously noted, the solution of the linear system was particularly fast and accurate with all test problems with closed boundary (Test Problems #1 through #4), or with open boundary but no inflow/outflow (Test Problem #5), while it was comparatively slower and less accurate with the test problems characterized by an open boundary and inflow/outflow (Test Problems #6 and #7). This is believed to be the main reason why the trends in Fig. 5 are stratified the way they are, with the results for Test Problems #6 and #7 somewhat removed from the cluster of the rest of the results. As previously noted, incompressible flow problems with inflow/outflow are intrinsically more difficult to handle at the discrete level than problems without inflow/outflow, and this should be borne in mind when making comparisons that include both types of problems.

Finally, convergence histories for ‖ $div(\underline{u}_h)$ ‖$_{L^2}$ and ‖ $div(\underline{u}_{hl})$ ‖$_{L^2}$ are provided in Fig. 6 for all test problems. As can be seen, the observed convergence rates are within 0.93-1.08 and appear consistent with $O(h)$ convergence. This is the same order of convergence of the velocity error in the H¹ norm, as is generally the case in numerical simulations [21].

## 6. Concluding remarks

The MINI mixed finite element discretization of the two-dimensional Stokes problem has been experimentally investigated on a selection of seven benchmark test cases with analytical solution, focusing on $O(h^{3/2})$ superconvergence and on the approximating properties of the complete computed velocity as compared to the piecewise-linear part of the computed velocity. Using automatically generated, unstructured triangular meshes, $O(h^{3/2})$ superconvergence of the pressure error $\| P - P_h \|_{L^2}$ and of the velocity error $\| i_h \underline{u} - \underline{u}_{hl} \|_{H^1}$ has been observed in all test problems, suggesting a validity of this type of superconvergence more general than covered by the existing theory, which is restricted to three-directional structured triangulations. Both the complete computed velocity $\underline{u}_h$ and the piecewise-linear part of the computed velocity $\underline{u}_{hl}$ approximated well the exact velocity $\underline{u}$ in all test problems. The complete computed velocity $\underline{u}_h$ is generally closer to the exact velocity (with the exception of flow problems with open boundary and inflow/outflow solved on coarse meshes, and only in the $H^1$ norm), and is therefore generally a better approximation than the piecewise-linear part of the computed velocity. The piecewise-linear part of the computed velocity $\underline{u}_{hl}$, on the other hand, conserves mass better and is therefore more appropriate for applications where the violation of the conservation of mass should be minimized. The present results will be instrumental in informing future studies on superconvergence and a-posteriori error estimation for automatic grid refinement, for the two-dimensional Stokes problem discretized with the MINI mixed finite element.

## Acknowledgement

The research of the second author (Daniele Boffi) has been partially supported by IMATI-CNR, GNCS-INDAM and has been performed within the MIUR initiative ''Dipartimenti di Eccellenza'' program (2018-2022).

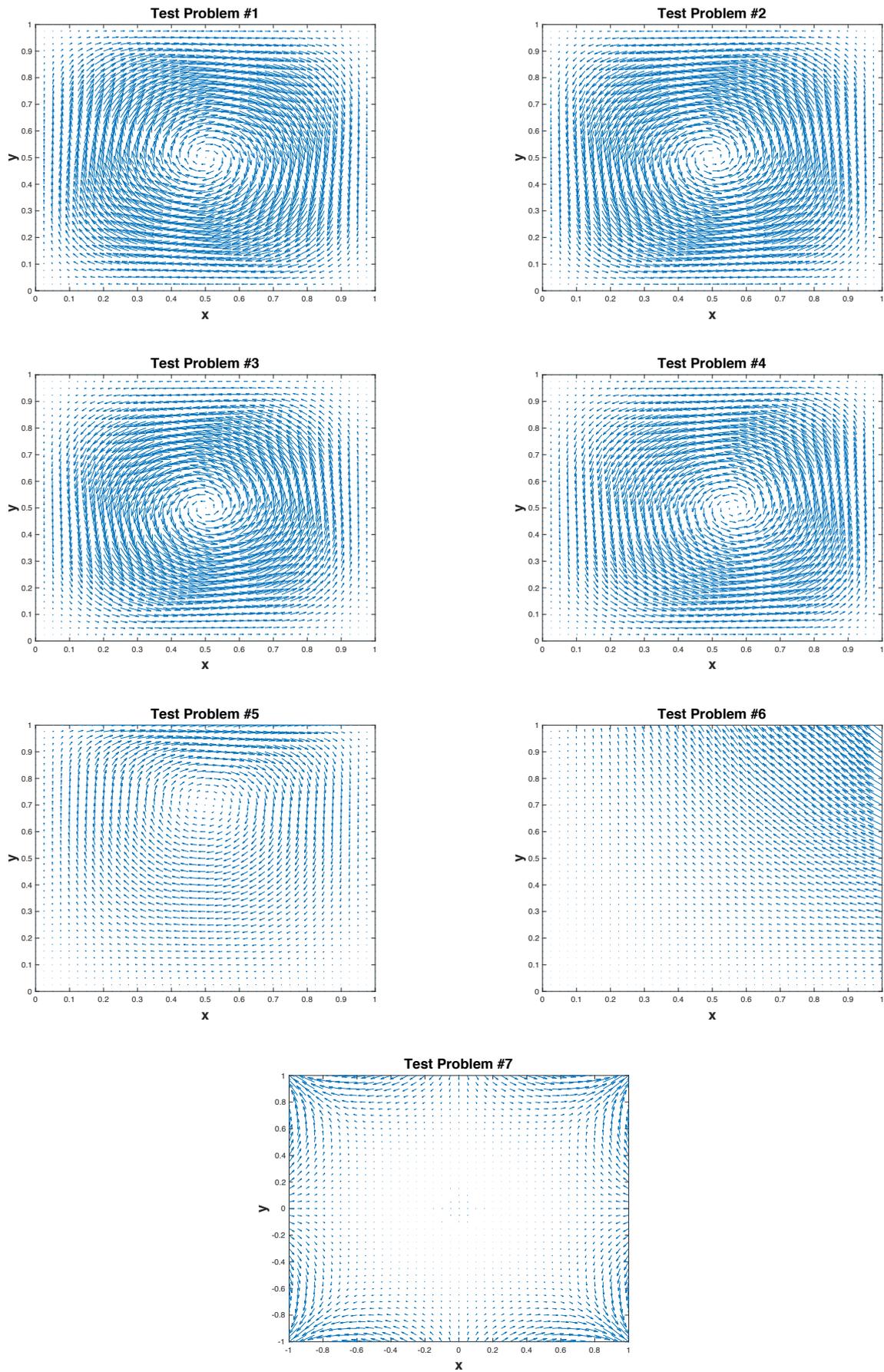

Fig. 1. Velocity vector fields for the benchmark test problems.

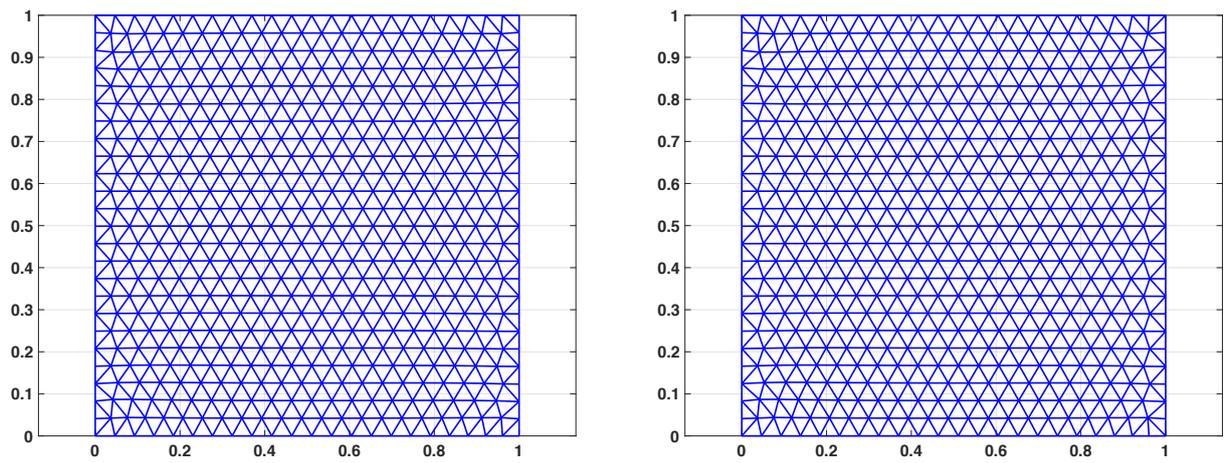

Fig. 2. Representative uniformly spaced mesh in the unit square (1014 triangles and 554 vertices) from DistMesh (left), an after correcting for corner triangles via diagonal exchange (right).

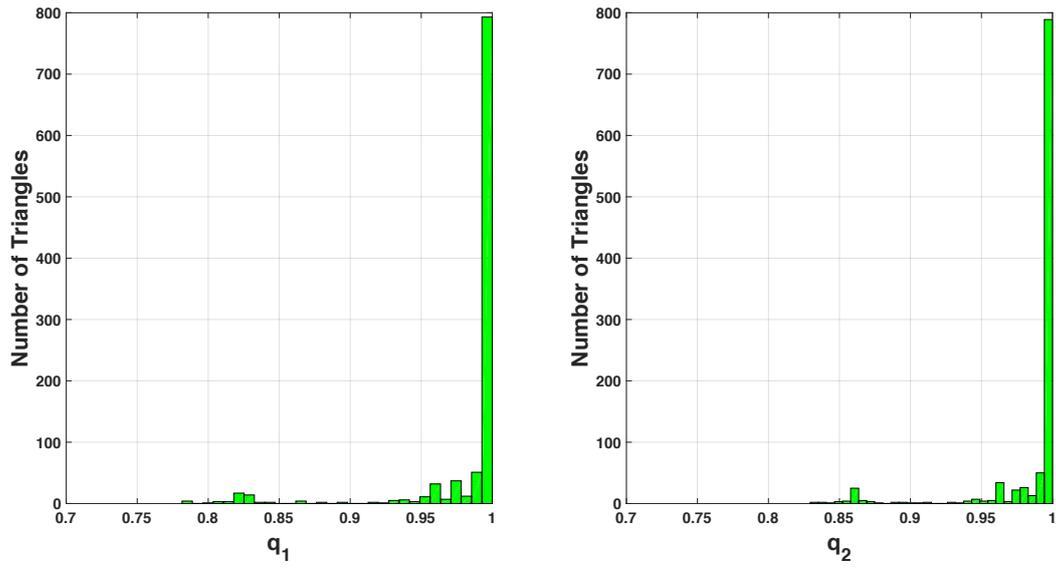

Fig. 3. Quality histograms for the mesh in Fig. 2 (right).

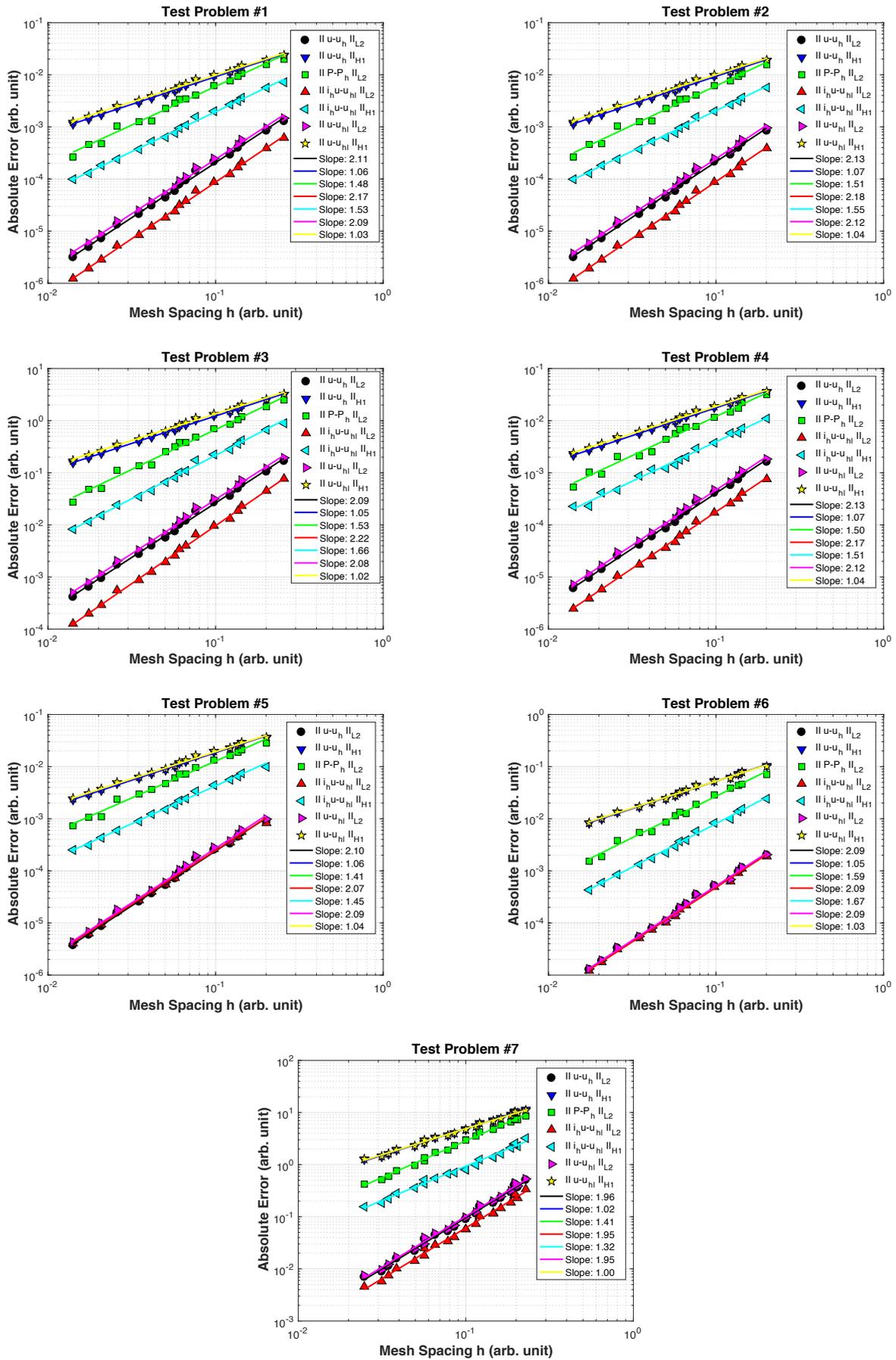

Fig. 4. Convergence histories for velocity and pressure errors (the solid lines are power law fits through the data points).

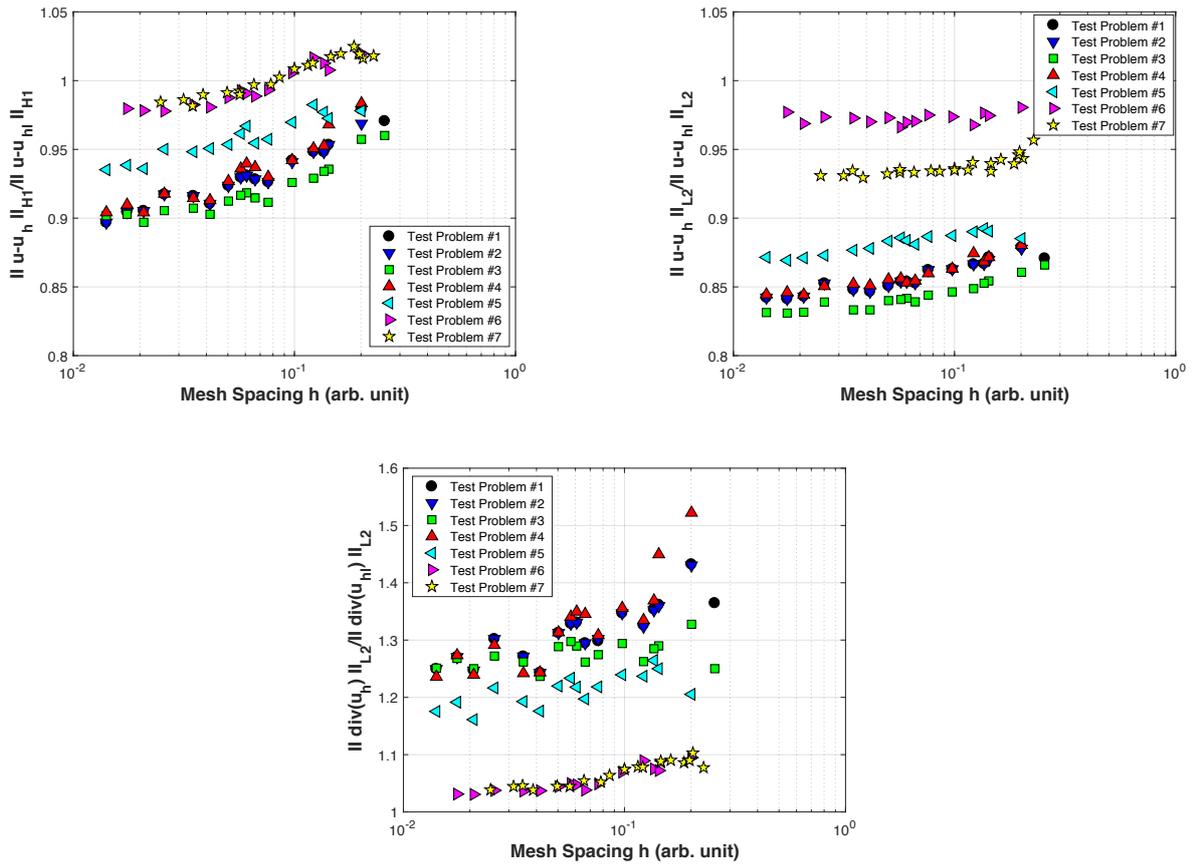

Fig. 5. Comparison between complete computed velocity $\underline{u}_h$ and piecewise-linear part of the computed velocity $\underline{u}_{hl}$. Top-left: ratio of velocity errors in $H^1$ norm; Top-right: ratio of velocity errors in $L^2$ norm; Bottom: ratio of $L^2$ norms of the divergence of the computed velocity.

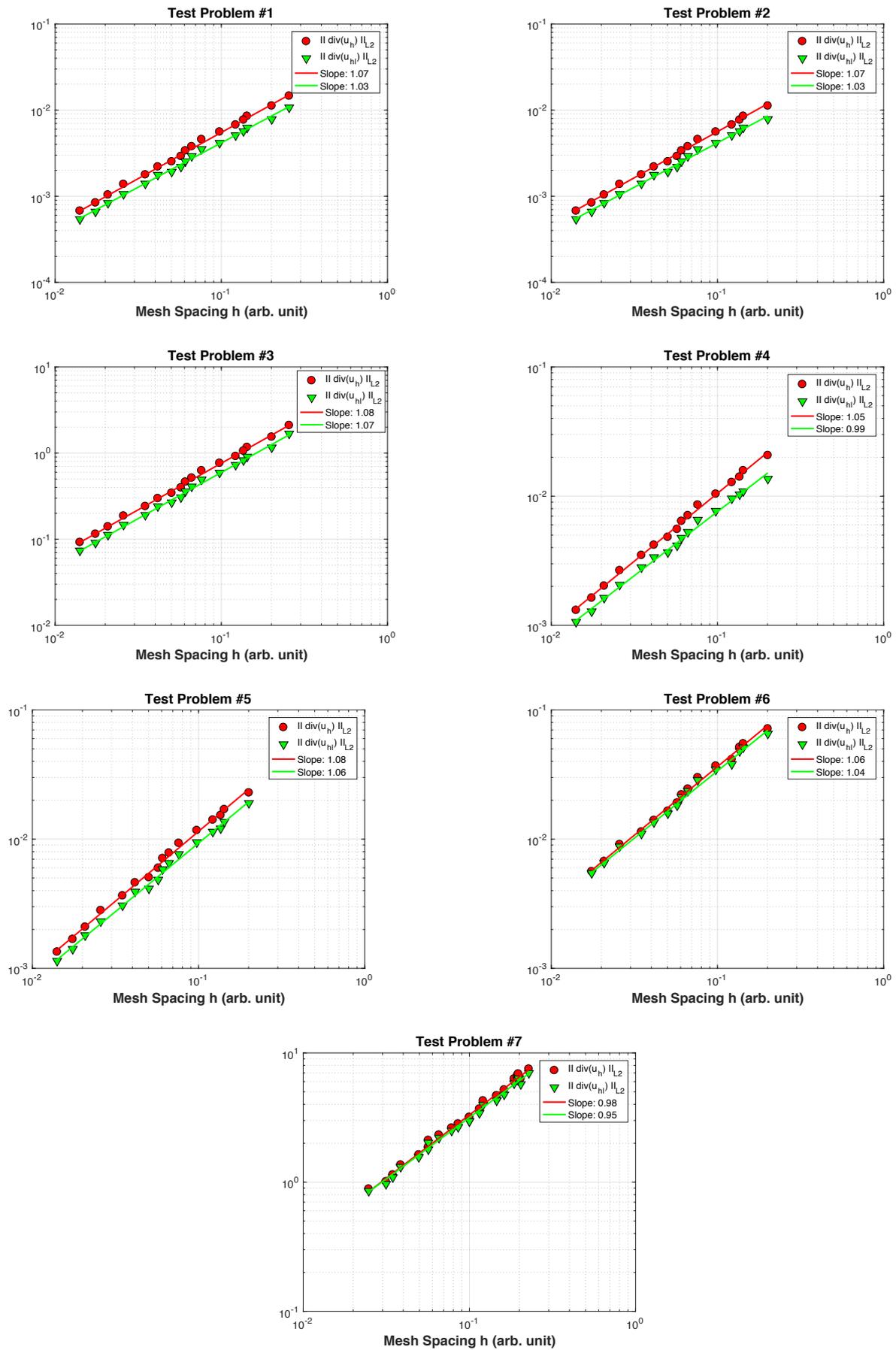

Fig. 6. Convergence histories for the $L^2$ norm of the divergence of the velocity (arbitrary unit on vertical axis; the solid lines are power law fits through the data points).

Table 1. Experimental Convergence Rates

| Problem # | $\|\underline{u} - \underline{u}_h\|_{L^2}$ | $\|\underline{u} - \underline{u}_h\|_{H^1}$ | $\|P - P_h\|_{L^2}$ | $\|i_h\underline{u} - \underline{u}_{hl}\|_{L^2}$ | $\|i_h\underline{u} - \underline{u}_{hl}\|_{H^1}$ | $\|\underline{u} - \underline{u}_{hl}\|_{L^2}$ | $\|\underline{u} - \underline{u}_{hl}\|_{H^1}$ |
|---|---|---|---|---|---|---|---|
| 1 | 2.11 | 1.06 | 1.48 | 2.17 | 1.53 | 2.09 | 1.03 |
| 2 | 2.13 | 1.07 | 1.51 | 2.18 | 1.55 | 2.12 | 1.04 |
| 3 | 2.09 | 1.05 | 1.53 | 2.22 | 1.66 | 2.08 | 1.02 |
| 4 | 2.13 | 1.07 | 1.50 | 2.17 | 1.51 | 2.12 | 1.04 |
| 5 | 2.10 | 1.06 | 1.41 | 2.07 | 1.45 | 2.09 | 1.04 |
| 6 | 2.09 | 1.05 | 1.59 | 2.09 | 1.67 | 2.09 | 1.03 |
| 7 | 1.96 | 1.02 | 1.41 | 1.95 | 1.32 | 1.95 | 1.00 |